\documentclass[a4paper,11pt]{article}
\usepackage{ams}
\usepackage[T1]{fontenc}

\newtheorem{df}{Definition}[section]
\newtheorem{thm}[df]{Theorem}
\newtheorem{prop}[df]{Proposition}
\newtheorem{cor}[df]{Corollary}
\newtheorem{ex}[df]{Example}
\newtheorem{lem}[df]{Lemma}

\def\boom{\quad\lower3pt\hbox{\vrule height1.1ex width .9ex depth -.2ex}
                    \vskip9pt}

\def\st{\ \vert\ }
\newcommand{\pf}{\noindent{\sc Proof.}\ }

\let\Bar=\overline

\let\Tilde=\widetilde

\def\tilalpha{\widetilde\alpha}
\def\tilbeta{\skew6\widetilde\beta}
\def\tiliota{\widetilde\iota}
\def\til0{\widetilde 0}
\def\tilone{\widetilde 1}
\def\tila{\widetilde a}

\def\tilq{\widetilde q}

\let\da=\partial
\let\isom=\cong
\let\sol=\bullet

\def\chigh{{\raise1.5pt\hbox{$\chi$}}}
\let\Ga=\Gamma
\let\Om=\Omega

\let\eps=\varepsilon
\let\phi=\varphi

\let\Ups=\Upsilon

\def\plusH{\ \lower 5pt\hbox{${\buildrel {\textstyle +}
\over {\scriptscriptstyle H}}$}\ }
\def\minusH{\ \lower 5pt\hbox{${\buildrel {\textstyle -}
\over {\scriptscriptstyle H}}$}\ }
\def\timesH{\ \lower 4pt\hbox{${\buildrel {\textstyle .}
\over{\scriptscriptstyle H}}$}\ }

\def\plusV{\ \lower 5pt\hbox{${\buildrel {\textstyle +}
\over {\scriptscriptstyle V}}$}\ }
\def\minusV{\ \lower 5pt\hbox{${\buildrel {\textstyle -}
\over {\scriptscriptstyle V}}$}\ }
\def\timesV{\ \lower 4pt\hbox{${\buildrel {\textstyle .}
\over{\scriptscriptstyle V}}$}\ }

\def\llangle{\langle\!\langle}
\def\rrangle{\rangle\!\rangle}

\def\LAgpd{${\cal LA}$-groupoid}
\def\LAdgpd{${\cal LA}$-double groupoid}
\def\PVBgpd{${\cal PVB}$-groupoid}
\def\VBgpd{${\cal VB}$-groupoid}

\def\id{{\rm id}}

\def\ddt#1{\left.\frac{d}{dt}#1\right|_0}

\def\pback#1{\mathbin{{{\lower1.2ex\hbox{$\times$}}\atop #1}}}
\def\sdp{\mathbin{\hbox{$\mapstochar\kern-.3333em\times$}}}

\def\gpd{\,\lower1pt\hbox{$\longrightarrow$}\hskip-.24in\raise2pt
             \hbox{$\longrightarrow$}\,}

\def\sgpd{\,\lower1pt\hbox{$\mlra$}\hskip-0.4in\raise2pt\hbox{$\mlra$}\,}
\def\svgpd{\Big\downarrow\!\!\Big\downarrow}

\def\vgpd{\Bigg\downarrow\!\!\Bigg\downarrow}

\def\vlra{\hbox{$\,-\!\!\!-\!\!\!-\!\!\!-\!\!\!-\!\!\!
-\!\!\!-\!\!\!-\!\!\!-\!\!\!-\!\!\!\longrightarrow\,$}}

\def\lgpd{\,\lower1pt\hbox{$\vlra$}\hskip-1.02in\raise2pt\hbox{$\vlra$}\,}

\def\lrah{\hbox{$\,-\!\!\!-\!\!\!
-\!\!\!-\!\!\!-\!\!\!-\!\!\!-\!\!\!\longrightarrow\,$}}

\def\mlra{\hbox{$\,-\!\!\!-\!\!\!\longrightarrow\,$}}

\def\hcompo#1#2#3{{\vcenter{\vbox{\hrule height.#2pt\hbox{\vrule width.#2pt
   height#1pt\kern#3pt\vrule width.#2pt\kern#3pt\vrule width.#2pt}
   \hrule height.#2pt}}}}

\def\hcomp{\mathchoice\hcompo643\hcompo643\hcompo{2.2}3{1.1}
                       \hcompo{1.4}3{.7}}

\def\vcompo#1#2#3{{\vcenter{\vbox{\hrule height.#2pt
   \hbox{\vrule width.#2pt height#3pt\kern#1pt\vrule width.#2pt}
   \hrule height.#2pt
   \hbox{\vrule width.#2pt height#3pt\kern#1pt\vrule width.#2pt}
   \hrule height.#2pt}}}}

\def\vcomp{\mathchoice\vcompo643\vcompo643\vcompo{2.2}3{1.1}
                       \vcompo{1.4}3{.7}}

\def\dsq{\mathop{\lower1pt\vbox{\hrule height.4pt \hbox
{\vrule width.4pt height.6em
\kern.6em \vrule width.4pt} \hrule height.4pt}}}

\def\dcomp{\mathop{\dsq\hskip-.88em\raise1pt\hbox{$\scriptstyle\nwarrow$}}}

\def\ssq{\vbox{\hrule height.4pt \hbox{\vrule width.4pt height.7in
\kern.7in \vrule width.4pt} \hrule height.4pt}}

\def\tsq{\mathop{\lower1pt\vbox{\hrule height.4pt \hbox
{\vrule width.4pt height.7em
\kern.7em \vrule width.4pt} \hrule height.4pt}}}

\def\m@th{\mathsurround=0pt}

\def\n@space{\nulldelimiterspace=0pt \m@th}

\def\Bbigg#1{{\hbox{$\left#1\vbox to 27.5pt{}\right.\n@space$}}}

\def\VGPD{\Bbigg\downarrow\!\!\Bbigg\downarrow}

\def\affel#1#2#3#4#5#6#7#8{\matrix{#7&\raise.7ex\hbox{$#1$}&#6\cr
                  \raise.5in\hbox{$#2$}&\sq&\raise.5in\hbox{$#3$}\cr
                          #8&#4&#5\cr}}

\begin{document}

\title{{\bf ON SYMPLECTIC DOUBLE GROUPOIDS AND THE DUALITY OF POISSON
GROUPOIDS}
\thanks{1991 {\em Mathematics
Subject Classification.} Primary 58H05. Secondary 17B66, 18D05, 22A22,
58F05.}}

\author{K. C. H. Mackenzie\\
        School of Mathematics and Statistics\\
        University of Sheffield\\
        Sheffield, S3 7RH\\
        England\\
        {\sf K.Mackenzie@sheffield.ac.uk}}

\date{{\sf May 19, 1998}}

\maketitle

\begin{abstract}
We prove that the cotangent of a double Lie groupoid $S$ has itself a
double groupoid structure with sides the duals of associated Lie
algebroids, and double base the dual of the Lie algebroid of the core
of $S$. Using this, we prove a result outlined by Weinstein in 1988, that the
side groupoids of a general symplectic double groupoid are Poisson
groupoids in duality. Further, we prove that any double Lie groupoid gives
rise to a pair of Poisson groupoids (and thus of Lie bialgebroids) in duality.
To handle the structures involved effectively we extend to this context
the dualities and canonical isomorphisms for tangent and cotangent
structures of the author and Ping Xu.
\end{abstract}

\newpage

In introducing the concept of Poisson groupoid, A.~Weinstein
\cite[\S4.5]{Weinstein:1988} defined two Poisson groupoids on the same base
manifold to be {\em dual} if the Lie algebroid dual of each is isomorphic
to the Lie algebroid of the other, and described a programme for showing
that, at least locally, Poisson groupoids in duality arise as the sides
of a symplectic double groupoid. In \cite{LuW:1989} J.--H. Lu and Weinstein
carried out this programme globally for Poisson groups, constructing for
a Poisson Lie group $G$ a symplectic double groupoid whose sides are $G$
and the dual group $G^*$. One purpose of the present paper is to explain for
the first time the general process, underlying \cite[\S4.5]{Weinstein:1988}
and \cite{LuW:1989}, by which a symplectic double groupoid gives rise to a
pair of Poisson groupoids in duality.

One would expect such properties of symplectic double groupoids to emerge as
a special case of the Lie theory for double groupoids proposed by the author.
In fact, whereas the construction of the double Lie algebroid of a double
Lie groupoid requires two steps \cite{Mackenzie:1992},
\cite{Mackenzie:Doubla2}, only the first is needed in order to understand
the duality associated with a symplectic double groupoid. The key is the
structure of the cotangent of a double Lie groupoid.

The fundamental example of an ordinary symplectic groupoid is the cotangent
groupoid $T^*G\gpd{\goth g}^*$ of a Lie group and, more generally,
the cotangent groupoid $T^*G\gpd A^*G$ of any Lie groupoid $G\gpd M$
\cite{CosteDW}. At the same time this groupoid, as a Lie groupoid, provides
a simple characterization of compatible Poisson structures on $G$: a Poisson
structure $\pi^\#\colon T^*G\to TG$ makes $G$ a Poisson groupoid if and only
if $\pi^\#$ is a morphism from $T^*G\gpd A^*G$ to $TG\gpd TM$
\cite{AlbertD:1991ssr}, \cite[8.1]{MackenzieX:1994}.

In \S\ref{sect:datcdg} of the present paper we define a corresponding
double structure on the cotangent of any double Lie groupoid $(S;H,V;M)$.
Here the side structures are $A^*_VS$ and $A^*_HS$, the duals of the Lie
algebroids of the two groupoid structures on $S$. Using the duality theory
for double structures of J.~Pradines \cite{Pradines:1988}, and the notion of
core for double groupoids introduced by R.~Brown and the author
\cite{BrownM:1992}, each of these is itself a groupoid over $A^*C$,
where $C$ is the core of $S$.

This result has two immediate consequences in \S\ref{sect:sdgapdg}.
It allows us to handle symplectic double groupoids ${\cal S}$,
and more generally Poisson double groupoids, in terms of the structure on
$T^*{\cal S}$: we prove (\ref{thm:sideduality}) by diagrammatic methods
that the side groupoids of a symplectic double groupoid are Poisson groupoids
in duality and that the core of the double groupoid provides a symplectic
realization of the double base. Secondly, for any double Lie groupoid $S$
we can take $T^*S$ itself as a symplectic double groupoid and it follows
that the Lie algebroid duals $A^*_VS$ and $A^*_HS$ are themselves Poisson
groupoids in duality (\ref{thm:pairs}). We will show in
\cite{Mackenzie:DLADLB} that this result is of crucial importance
in understanding general double Lie algebroids.

In \S\ref{sect:dflag} we are ultimately concerned to extend the known
duality between Lie algebroids and Poisson vector bundles
\cite{Courant:1990} to \LAgpd s and what we call \PVBgpd s.
This unexpectedly difficult result (\ref{thm:LAPVB}) is
needed in \ref{thm:needed}. To work with a Poisson structure on a \VBgpd\
requires the cotangent of the \VBgpd\ and this has a triple structure which
is not straightforward. To handle this we extend to general \VBgpd s the
dualities and canonical isomorphisms associated with the tangent and
cotangent structures of a Lie groupoid which were introduced by the author
and Ping Xu in \cite{MackenzieX:1994} and \cite{MackenzieX:1998}. In the
process we obtain (\ref{prop:DRI}) a concrete description of the duality
between the two bialgebroids of a double Lie groupoid; this will be important
in \cite{Mackenzie:DLADLB}.

An early version of the material of this paper was presented at the
Groupoidfest at Berkeley in 1997. I am grateful to Alan Weinstein and
Ping Xu for comments there, and to Johannes Huebschmann for conversations
at the December 1997 Workshop at Warwick.

\newpage

\section{THE COTANGENT DOUBLE GROUPOID}
\label{sect:datcdg}

We begin by recalling the notion of dual for \VBgpd s due to Pradines
\cite{Pradines:1988}. A \VBgpd\ $(\Om;G,A;M)$ (see also
\cite[\S4]{Mackenzie:1992}) is a double Lie groupoid in which two opposite
structures are vector bundles;
\begin{equation}                                      \label{eq:VBgpd}
\matrix{                    &&\tilq &&\cr
                            &\Om&\lrah &G&\cr
                            &&&&\cr
        \tilalpha, \tilbeta &\vgpd&&\vgpd&\alpha,\beta\cr
                            &&&&\cr
                            &A  &\lrah &M,&\cr
                            &&q&&\cr}
\end{equation}
that is, $\Om$ is a vector bundle over $G$,
which is a Lie groupoid over $M$, and $\Om$ is also a Lie groupoid over $A$,
which is a vector bundle over $M$, subject to the condition that the structure
maps of the groupoid structure (source, target, identity, multiplication,
inversion) are vector bundle morphisms, and the \lq\lq double source
map\rq\rq\ $(\tilq,\tilalpha)\colon\Om\to G\pback{M}A$ is a surjective
submersion.

The core
$K = \{\xi\in\Om\st \tilalpha(\xi) = 0_m,\ \tilq(\xi) = 1_m,\ \exists m\in M\}$
of $\Om$ is then a vector bundle over $M$ under the restriction of the
vector bundle operations on $\Om$. (Here one can also regard $\Om$ as an
\LAgpd, as in \cite[\S5]{Mackenzie:1992}, in which both Lie algebroid
structures are zero.)

Let $\Om^*$ be the dual of $\Om$ as a vector bundle over $G$. Define a
groupoid structure on $\Om^*$ with base $K^*$ as follows. Take $\Phi\in
\Om^*_g$ where $g\in G_m^n$. Then the source and target of $\Phi$ in
$K^*_m$ and $K^*_n$ respectively are
\begin{equation}                                    \label{eq:dual1}
\langle\tilalpha_*(\Phi),k \rangle =
         \langle\Phi,-\til0_gk^{-1}\rangle, \quad k \in K_m,
\qquad\qquad
\langle\tilbeta_*(\Phi),k \rangle =
         \langle\Phi, k \til0_g\rangle, \quad k \in K_n.
\end{equation}
For the composition, take $\Psi\in\Om^*_h$ with
$\tilalpha_*(\Psi) = \tilbeta_*(\Phi)$. Any element of $\Om_{hg}$ can
be written as a product $\eta\xi$ where $\eta\in\Om_h$ and $\xi\in\Om_g$.
Now the compatibility condition on $\Psi$ and $\Phi$ ensures that
\begin{equation}                                    \label{eq:dual2}
\langle\Psi\Phi, \eta\xi\rangle = \langle\Psi,\eta\rangle +
                                    \langle\Phi,\xi\rangle
\end{equation}
is well defined. The identity element of $\Om^*$ at $\theta\in K^*_m$ is
$\tilone_\theta\in\Om^*_{1_m}$ defined by
\begin{equation}                                    \label{eq:dual3}
\langle\tilone_\theta,\tilone_X + k \rangle = \langle\theta,k \rangle,
\end{equation}
where any element of $\Om_{1_m}$ can be written as $\tilone_X + k $
for some $X\in A_m$ and $k \in K_m$.

It is straightforward to check that $(\Om^*;G,K^*;M)$ is a \VBgpd, the
{\em dual \VBgpd\ to $\Om$}. The core of $\Om^*$ is the vector bundle
$A^*\to M$, with the core element corresponding to $\phi\in A^*_m$ being
$\Bar\phi\in\Om^*_{1_m}$ defined by
\begin{equation}                                       \label{eq:dual4}
\langle\Bar\phi,\tilone_X + k \rangle =
   \langle\phi,X + \da_A(k )\rangle
\end{equation}
for $X\in A_m,\ k \in K_m$. Here $\da_A\colon K\to A$ is the restriction of
$\tilbeta\colon\Om\to A$; it is a vector bundle morphism over $M$ (see
\cite[\S5]{Mackenzie:1992}).

Note that the dual of $\Om^*$ identifies canonically with $\Om$ as a
\VBgpd, the signs in (\ref{eq:dual1}) cancelling. The $\da$ map $A^*\to K^*$
for $\Om^*$ is $\da_A^*$.

Now consider a morphism of \VBgpd s which preserves the lower groupoids:
$$
(F;\id_G,f;\id_M)\colon(\Om;G,A,M)\to (\Om';G,A';M)
$$
and denote the core morphism $K\to K'$ by $f_K$. The proof of the following
result is simple.

\begin{prop}                                  \label{prop:dualmorph}
The dual morphism $F^*\colon\Om'^*\to\Om^*$ is a morphism of the dual
\VBgpd s, with base map $f_K^*\colon K'^*\to K^*$ and core morphism
$f^*\colon A'^*\to A^*$.
\end{prop}

The two basic examples are the tangent and the cotangent \VBgpd s of a
Lie groupoid. For any Lie groupoid $G\gpd M$, the tangent groupoid has a
\VBgpd\ structure $(TG;G,TM;M)$ with core $AG\to M$. The map
$\da\colon AG\to TM$ is the anchor $a$ of $AG$.

Taking the dual gives the {\em cotangent groupoid} $(T^*G;G,A^*G;M)$ of
\cite{CosteDW} or \cite{Pradines:1988}. Here the source and target maps are
given for $\Phi\in T^*_gG$ by
\begin{equation}                                    \label{eq:T^*G}
\Tilde{\alpha}(\Phi)(X) = \Phi(T(L_g)(X - T(1)(a(X))),\qquad
\Tilde{\beta}(\Phi)(Y) = \Phi(T(R_g)(Y)),
\end{equation}
where $X\in A_{\alpha g}G,\ Y\in A_{\beta g}G$; thus
$\Tilde{\alpha}(\Phi)\in A^*_{\alpha g}G$ and
$\Tilde{\beta}(\Phi)\in A^*_{\beta g}G$.
Here $L_g$ and $R_g$ are the left and right translations in $G$. If $\Psi
\in T^*_hG$ and $\Tilde{\alpha}(\Psi) = \Tilde{\beta}(\Phi)$ then
$\alpha h = \beta g$ and we define $\Psi\Phi\in T^*_{hg}G$ by
\begin{equation}                                    \label{eq:dualcomp}
(\Psi\Phi)(Y\sol X) = \Psi(Y) + \Phi(X),\qquad Y\in T_hG,\ X\in T_gG.
\end{equation}
If $\phi\in A^*_mG$, then the identity element over $\phi$ is $\Tilde{1}_\phi
\in T^*_{1_m}G$ defined by $\Tilde{1}_\phi(T(1)(x) + X) = \phi(X)$ for $X\in
A_mG,\ x\in T_m(M)$. The core is $T^*M$, where
$\omega\in T^*_m(M)$ defines the core element
$\Bar{\omega}\in T^*_{1_m}G$ given by $\Bar{\omega}(T(1)(x) + \Bar{X})
= \omega(x + a(X))$ for $x\in T_m(M),\ X\in A_mG$. Lastly, the map
$\da\colon T^*M\to A^*G$ is $a^*$. This is precisely the cotangent
structure as used in \cite[\S7]{MackenzieX:1994}.

\begin{ex}\rm                                  \label{ex:rmk31}
If now $G\gpd M$ is a Poisson groupoid, the Poisson tensor
$\pi_G^\#\colon T^*G\to TG$ is a morphism of \VBgpd s with respect to the
base maps $\id_G$ and $a_*\colon A^*G\to TM$ (see \cite{AlbertD:1991ssr} or
\cite[8.1]{MackenzieX:1994}). It now follows from \ref{prop:dualmorph} and
from the skewsymmetry of $\pi_G^\#$, that the core map of $\pi_G^\#$ is
$-a_*^*$.
\end{ex}

We turn now to the double cotangent groupoid. Consider a double Lie
groupoid $(S;H,V;M)$ as defined in \cite[\S2]{Mackenzie:1992}. Thus $S$ has
two Lie groupoid structures, a {\em horizontal structure} with base $V$,
and a {\em vertical structure} with base $H$, where both $V$ and $H$ are
Lie groupoids on base $M$, such that the structure maps of each groupoid
structure on $S$ are morphisms with respect to the other. We display
$(S;H,V;M)$ and a typical element as in Figure~\ref{fig:S}. It is further
assumed that the double source map
$(\tilalpha_V,\tilalpha_H)\colon S\to H\pback{M}V$ is a surjective
submersion.

\medskip

\begin{figure}[htb]
\begin{picture}(395,100)
\put(30,50){$\matrix{&&{\tilalpha_H,\tilbeta_H}&&\cr
          &S&\lgpd &V&\cr
          &&&&\cr
          {\tilalpha_V,\tilbeta_V}&\vgpd&&\vgpd&{\alpha_V,\beta_V}\cr
          &&&&\cr
          &H&\lgpd &M&\cr
          &&{\alpha_H,\beta_H}&&\cr}$}
\put(300,50){$\matrix{&\raise.7ex\hbox{$\tilbeta_V(s)$}&\cr
                  \raise.35in\hbox{$\tilbeta_H(s)$}&\ssq&\raise.35in
                         \hbox{$\tilalpha_H(s)$}\cr
                         &\tilalpha_V(s)&\cr}$}
\put(365,55){$s$}
\end{picture}         \caption{\ \label{fig:S}}
\end{figure}

The core $C$ of $S$ is the manifold of elements
$c\in S$ with $\tilalpha_H(c) = \tilone^V_m,\ \tilalpha_V(c)= \tilone^H_m$
for some $m\in M$. The core has a groupoid structure over $M$ with source
$\alpha_C(c) = \alpha_V(\tilalpha_H(c))$, target
$\beta_C(c) = \beta_V(\tilbeta_H(c))$, and composition
$$
c' \dcomp c = (c' \ \hcomp\ \tilone^V_{h_2})\ \vcomp\ c =
                     (c' \ \vcomp\ \tilone^H_{v_2})\ \hcomp\ c
$$
where $v_2 = \tilbeta_H(c),\ h_2 = \tilbeta_V(c).$ The identity of $C$ at
$m \in M$ is $1^C_m = 1^2_m$ and the inverse of $c \in C$ is
\begin{equation}                              \label{eq:Cinv}
c^{-1(C)} = c^{-1(H)}\ \vcomp\ \tilone^H_{v^{-1}}
            = c^{-1(V)}\ \hcomp\ \tilone^V_{h^{-1}}\,,
\end{equation}
where $v = \tilbeta_H(c)$ and $h = \tilbeta_V(c)$. With this structure the
maps $\da_H\colon C\to H,\ c\mapsto\tilbeta_V(c)$ and
$\da_V\colon C\to V,\ c\mapsto\tilbeta_H(c)$ are morphisms over $M$. See
\cite[\S2]{BrownM:1992} or \cite[\S2, p.197]{Mackenzie:1992}.

Applying the Lie functor to the vertical structure $S\gpd H$ yields a Lie
algebroid $A_VS\to H$; because the Lie functor preserves pullbacks, the
horizontal groupoid $S\gpd V$ prolongs to a groupoid structure
$A_VS\gpd AV$. These two structures, shown in Figure~\ref{fig:LAgpd}(a),
together constitute what we called in \cite[\S4]{Mackenzie:1992} an \LAgpd.
Similarly we can take the horizontal Lie algebroid $A_HS\to V$; it is a
groupoid over $AH$ as in Figure~\ref{fig:LAgpd}(b).
\begin{figure}[htb]
\begin{picture}(395,130)
\put(0,70){$\matrix{&&A(\tilalpha_H),A(\tilbeta_H)&&\cr
        &A_VS  &\lgpd&AV&\cr
        &&&&\cr
\tilq_V &\Bigg\downarrow&&\Bigg\downarrow&\cr
        &&&&\cr
        &H  &\lgpd&M&\cr
        &&\alpha_H,\beta_H&&\cr}$}
\put(70,70){$TH\stackrel{T(\alpha_H),T(\beta_H)}{\lgpd}\ TM$}
\put(45,95){\vector(1,-1){15}}    \put(60,85){$\tila_V$}
\put(160,95){\vector(1,-1){15}}           
\put(65,65){\vector(-1,-1){15}}
\put(175,65){\vector(-1,-1){15}}
\put(70,0){(a)}
\put(260,70){$\matrix{      &&\tilq_H &&\cr
                            &A_HS    &\vlra &V &\cr
                            &&&&\cr
                            &\vgpd&&\vgpd&\cr
                            &&&&\cr
                            &AH &\vlra &M&\cr
                            &&&&\cr}$}
\put(340,35){$\VGPD$}
\put(290,95){\vector(2,-1){35}}    \put(300,70){$\tila_H$}
\put(290,30){\vector(2,-1){35}}    \put(300,10){$a_H$}
\put(355,80){\vector(2,1){35}}     \put(335,65){$TV$}
\put(360,5){\vector(2,1){35}}      \put(335,0){$TM$}
   \put(350,-20){(b)}
\end{picture}         \caption{\ \label{fig:LAgpd}}
\end{figure}
An {\em \LAgpd}\ is a \VBgpd\ as in (\ref{eq:VBgpd}) in which both $\Om\to G$
and $A\to M$ have Lie algebroid structures and each of the groupoid structure
maps is a Lie algebroid morphism. Each $k\in\Ga K$ induces a section
$\Bar{k}\in\Ga_G\Om$ where $\Bar{k}(g) = k(\beta g)\til0_g$; define a
bracket on $\Ga K$ by $[\Bar{k}, \Bar{\ell}] = \Bar{[k,\ell]}.$ Then $K$ is
a Lie algebroid on $M$ with anchor $a_K = a\circ\da_A$, and both $\da_A$ and
the restriction $\da_{AG}$ of the anchor $\tila\colon\Om\to TG$ to $K\to AG$,
are Lie algebroid morphisms. In particular $a_K = a_{AG}\circ\da_{AG}$
also. See \cite[\S5]{Mackenzie:1992}.

We quote the following result from \cite[1.6]{Mackenzie:Doubla2}.

\begin{thm}                               \label{thm:cores}
Let $(S;H,V;M)$ be a double Lie groupoid with core groupoid $C$. Then
$AC$, the Lie algebroid of $C$, is canonically isomorphic to both the
core Lie algebroid of the \LAgpd\ $A_VS$ and the core Lie algebroid of
the \LAgpd\ $A_HS$.
\end{thm}

The two groupoid structures on $S$ give rise to two cotangent groupoid
structures on $T^*S$, one over $A^*_HS$ and one over $A^*_VS$. In turn,
since both $A_HS$ and $A_VS$ are \VBgpd s with core $AC$, the duals
$A^*_HS$ and $A^*_VS$ both have groupoid structures on base $A^*C$.

\begin{thm}                                         \label{thm:dcot}
With the structures just described, $T^*S$ is a double Lie groupoid
\begin{equation}                       \label{eq:doublecot}
\matrix{&&{\tilalpha_{*H},\tilbeta_{*H}}&&\cr
          &T^*S&\lgpd &A^*_HS&\cr
          &&&&\cr
          {\tilalpha_{*V},\tilbeta_{*V}}&\vgpd&&\vgpd&{\alpha_{*V},\beta_{*V}}\cr
          &&&&\cr
          &A^*_VS&\lgpd &A^*C&\cr
          &&{\alpha_{*H},\beta_{*H}}&&\cr}
\end{equation}
\end{thm}

We omit the proof of this, in favour of the proof of the following theorem.
The proof of \ref{thm:dcot} is a long but straightforward verification,
using the groupoid structures defined in (\ref{eq:dual1})---(\ref{eq:dual3})
on the side groupoids and, on $T^*S$, the cotangent groupoid
structures (\ref{eq:T^*G}).

The following theorem is somewhat surprising, in that the cotangent behaves
as if it were covariant.

\begin{thm}
The core groupoid of {\em (\ref{eq:doublecot})} is naturally isomorphic
with the cotangent groupoid $T^*C\gpd A^*C$ of the core of $S$.
\end{thm}

\pf
Define a map $E\colon T_CS\to TC$ which sends each vector tangent to $S$
at a point of the core, to a vector at the same point tangent to the core
itself. Take $c\in C$ with $v = \da_V(c),\ h = \da_H(c)$, where $v\in V_m,\
h\in H_m$. Represent an element $\xi\in T_cS$ as
$$
\matrix{w&\raise.7ex\hbox{$W$}&x\cr
                  \raise.35in\hbox{$Z$}&\ssq&\raise.35in\hbox{$X$}\cr
                  z&Y&y\cr}\hskip-.7in\raise5pt\hbox{$\xi$}
$$
where $X = T(\tilalpha_H)(\xi)\in T_{1^V_m}V,\ x = T(\beta_V)(X) =
T(\alpha_H)(W) \in T_mM$, et cetera, and define
\begin{equation}
E(\xi) = \xi - T(L^V_c)T(\tilone^H)(X - T(1^V)(x))
             - T(L^H_c)T(\tilone^V)(Y - T(1^H)(z)),
\end{equation}
where $L^V_c$ and $L^H_c$ denote left translation in the two groupoid
structures on $S$.

Take $\sigma\in T^*_cC$ and define $\Sigma\in T^*_cS$ by
$\Sigma(\xi) = \sigma(E(\xi)),\ \xi\in T_cS$. We must show that $\Sigma$
is a core element. Denote $\tilalpha_{*V}(\Sigma)\in A^*_VS|_{1^H_m}$ by
$\rho$, and $\alpha_{*H}(\rho)\in A^*_mC$ by $\theta$. Then for all
$k \in A_mC$,
$$
\langle\theta,k \rangle = -\rho(T(\tiliota_H)(k ))
= \Sigma(T(L^V_c)T(\tiliota_V)T(\tiliota_H)(k ))
$$
where $\tiliota_H$ and $\tiliota_V$ are the two inversions in $S$.

We first prove that $\rho = 1^{*V}_\theta$. Take any
$\xi = A(\tilone^H)(X) + k \in A_VS|_{1^H_m}$ where $X\in A_mV,\
k \in A_mC$. Then $1^{*V}_\theta(\xi) = \langle\theta,k \rangle$
and $\rho(\xi) = - \Sigma(T(L^V_c)T(\tiliota_V)(\xi))$. Now although it is
natural to consider $k $ as an element of $A_VS$ in this context, it may
also be regarded as an element of $A_HS$ and we then get
$$
T(\tiliota_H)(k ) = T(\tilone^H)(Z) - k
$$
where $Z = \da_{AV}(k )$. Next, applying $E$ to
$T(L^V_c)T(\tiliota_V)T(\tilone^H)(Z)$, or equally with $X$ in place of $Z$,
gives $0\in T_cC$, where we use
$\tiliota_V\circ\tilone^H = \tilone^H\circ\iota_V$. Putting these together,
the equation for $\rho$ follows.

We must also prove that $\theta\in A^*_mC$ is the source of $\sigma$ with
respect to $T^*C\gpd A^*C$. For $k \in A_mC$ the above shows that
$$
\langle\theta,k \rangle = -\Sigma(T(L^V_c)T(\tiliota_V)(k )).
$$
Calculating $E$ of the argument on the right hand side, we obtain
\begin{equation}                             \label{eq:sum}
T(L^V_c)T(\tiliota_V)(k ) + T(L^H_c)T(\tilone^V)(T(1^H)(a_HW) - W),
\end{equation}
where $W = \da_{AH}(k )$. Write $k  = \ddt{c_t}$ where the curve
$c_t$ in $C$ is of the form
$$
\matrix{ &\raise.7ex\hbox{$h_t$}&\cr
                  \raise.35in\hbox{$v_t$}&\ssq&\raise.35in\hbox{$1^V_m$}\cr
                   &1^H_m&\cr}\hskip-.7in\raise5pt\hbox{$c_t$}
$$
Then the first term in (\ref{eq:sum}) is $\ddt{c\,\vcomp\, c_t^{-1(V)}}$ and
the second is $\ddt{c\,\hcomp\, \tilone^V_{h_t^{-1}}}$. The sum is therefore
the derivative of the product shown in Figure~\ref{fig:squares},
\begin{figure}[htb]
\begin{picture}(350,120)  
\put(160,65){$\ssq$}                    \put(230,65){$\ssq$}
\put(180,85){$c$}                       \put(255,85){$\tilone^V_{h_t^{-1}}$}
\put(160,0){$\ssq$}                     \put(230,0){$\ssq$}
\put(180,25){$c_t^{-1(V)}$}             \put(255,25){$\tilone^V_{h_t^{-1}}$}
\end{picture}\caption{\ \label{fig:squares}}
\end{figure}
and in terms of the multiplication $\dcomp$ and the inversion
$c\mapsto c^{-1(C)}$ in $C\gpd M$, this product is $c\dcomp c_t^{-1(C)}$
(see (\ref{eq:Cinv})). The
derivative is therefore $T(L^C_c)T(\iota_C)(k )$ and so we have
$\langle\theta,k \rangle = \langle\alpha_{*C}(\sigma),k \rangle$.

The proof that $\tilalpha_{*H}(\Sigma) = 1^{*H}_\theta$ now follows in
the same way. One likewise checks that $\beta_{*C}(\sigma) =
\beta_{*H}(\tilbeta_{*V}(\Sigma))$ and that the multiplications
correspond.

This proves that $\sigma\mapsto\Sigma$ is an isomorphism into the core of
$T^*S$. We leave the reader to check that the image is the whole of $T^*C$.
\boom

In point of fact, $T^*S$ is a triple structure, each of its four spaces being
also a vector bundle over the corresponding space of $(S;H,V;M)$, and all the
groupoid structure maps being vector bundle morphisms. We will return to this
in \S\ref{sect:sdgapdg}.

\section{SYMPLECTIC DOUBLE GROUPOIDS AND POISSON DOUBLE GROUPOIDS}
\label{sect:sdgapdg}

As with ordinary groupoids, it is actually easier to study symplectic
structures by specializing from the Poisson case.
In what follows we will repeatedly use the following simple result.

\begin{prop}
\begin{enumerate}
\item Let $(\phi;\phi_H,\phi_V;\phi_M)\colon (S;H,V;M)\to (S';H',V';M')$
be a morphism of double groupoids with core morphism $\phi_C\colon C\to C'$.
Then $\da_H'\circ\phi_C = \phi_H\circ\da_H$ and $\da_V'\circ\phi_C =
\phi_V\circ\da_V$.

\item Let $(\phi;\phi_G,\phi_A;\phi_M)\colon (\Om;G,A;M)\to (\Om';G',A';M')$
be a morphism of \LAgpd s with core Lie algebroid morphism
$\phi_K\colon K\to K'$. Then $\da_{AG}'\circ\phi_K = A(\phi_G)\circ\da_{AG}$
and $\da_A'\circ\phi_K = \phi_A\circ\da_A$.
\end{enumerate}
\end{prop}

\begin{df}
A {\em Poisson double groupoid} is a double Lie groupoid
$({\cal S};{\cal H},{\cal V};P)$ together with a Poisson structure on
${\cal S}$ with respect to which both groupoid structures on ${\cal S}$ are
Poisson groupoids.
\end{df}

For the theory of ordinary Poisson groupoids, see \cite{Weinstein:1988},
\cite{MackenzieX:1994},  \cite{Xu:1995}, \cite{MackenzieX:1998} and, for
Lie bialgebroids, \cite{Kosmann-Schwarzbach:1995}. For a Poisson groupoid
$G\gpd P$ we take the Poisson structure on the base to be
$\pi_P^\# = a_*\circ a^* = -a\circ a_*^*$. This convention (used in
\cite{Kosmann-Schwarzbach:1995}) is opposite to that used in
\cite{MackenzieX:1994}, but is necessary in order that the two $\da$ maps
$a^*\colon T^*P\to A^*G$ and $-a^*_*\colon T^*P\to AG$ of the cotangent
\LAgpd\ be Lie algebroid morphisms: see \ref{ex:rmk31}.

The two Poisson groupoid structures on ${\cal S}$ induce maps
$\tila_{*{\cal H}}\colon A^*_V{\cal S}\to T{\cal H}$ and
$\tila_{*{\cal V}}\colon A^*_H{\cal S}\to T{\cal V}$, the
anchors for the Lie algebroid structures on the duals of
$A_V{\cal S}\to {\cal H}$ and $A_H{\cal S}\to {\cal V}$,
with respect to which $\pi_{\cal S}^\#\colon T^*{\cal S}\to T{\cal S}$ is a
morphism of each of the two groupoid structures on $T^*{\cal S}$. From the
following result it follows that $\pi^\#_{\cal S}$ is actually a morphism
of double groupoids over a map $a_{*{{\cal C}}}\colon A^*{\cal C}\to TP$.

\begin{lem}                                      \label{lem:simple}
Let $(S;H,V;M)$ and $(S';H',V';M')$ be double Lie groupoids and let
$\phi\colon S\to S',\ \phi_H\colon H\to H'$ and $\phi_V\colon V\to V'$ be
maps such that $(\phi,\phi_H)$ and $(\phi,\phi_V)$ are morphisms of the
two ordinary groupoid structures on $S$ and $S'$. Then there is a unique
map $\phi_M\colon M\to M'$ such that $(\phi;\phi_H,\phi_V;\phi_M)$ is a
morphism of double groupoids.
\end{lem}

\pf
Take $m\in M$. The double identity $1^2_m$ can be written both as
$\tilone^V_{1^H_m}$ and as $\tilone^H_{1^V_m}$. Its image under $\phi$
is therefore an identity for both top structures on $S'$, and must
therefore be a double identity $1^2_{\phi_M(m)}$. It also follows that
$\phi_H,\ \phi_V$ and $\phi_M$ commute with the source and target
projections.

Now since $\phi(\tilone^V_h) = \tilone^V_{\phi_H(h)}$ for all $h\in H$
and $\tilone^V_{h_1}\,\hcomp\,\tilone^V_{h_2} = \tilone^V_{h_1h_2}$
for all compatible $h_1,h_2\in H$, it follows that $\phi_H$ is a morphism
over $\phi_M$. Similarly for $\phi_V$.
\boom

Returning to the Poisson double groupoid ${\cal S}$, the bases ${\cal H}$
and ${\cal V}$ acquire Poisson structures which we take to be
$\pi_{\cal H}^\# = \tila_{*{{\cal H}}}\circ\tila^*_{{\cal V}}$ and
$\pi_{\cal V}^\# = \tila_{*{{\cal V}}}\circ\tila^*_{{\cal H}}$. We will
prove below that ${\cal H}$ and ${\cal V}$ are Poisson groupoids with
respect to these structures.

First note that the core morphism $T^*{\cal C}\to T{\cal C}$ of
$\pi_{\cal S}^\#$ defines a Poisson structure on ${\cal C}$; denote this
by $\pi_{\cal C}^\#$. The core morphism is a morphism of groupoids over
$a_{*{{\cal C}}}\colon A^*{\cal C}\to TP$ and it follows that ${\cal C}\gpd P$
is a Poisson groupoid. Give $P$ the Poisson structure
$\pi_P^\# = a_{*{\cal C}}\circ a^*_{{\cal C}}$ induced from ${\cal C}$.

Since ${\cal S}_V\gpd {\cal H}$ is a Poisson groupoid, its Lie
algebroid dual $A_V^*{\cal S}$ has a Lie algebroid structure with
anchor $\tila_{*{{\cal H}}}$. Similarly $A^*{\cal C}$ has a Lie algebroid
structure with anchor $a_{*{\cal C}}$. The next result will be proved at
the end of \S\ref{sect:dflag}.

\begin{thm}                                 \label{thm:needed}
With respect to these structures
$(A^*_V{\cal S};{\cal H},A^*{\cal C};P)$ is an \LAgpd.
\end{thm}

We may therefore consider both $T^*{\cal S}$ and $T{\cal S}$ to be triple
structures as in Figure~\ref{fig:triples}: to be precise, they are double
groupoid objects in the category of Lie algebroids.
\begin{figure}[htb]
\begin{picture}(350,200)  
\put(0,150){$\matrix{&&      &\cr
                      &T^*{\cal S}&\sgpd &A^*_H{\cal S}\cr
                      &&&\cr
                      &\svgpd          & &\svgpd          \cr
                      &&&\cr
                      &A^*_V{\cal S}&\sgpd &A^*{\cal C}\cr}$}

\put(40, 160){\vector(3,-4){40}}                
\put(110, 160){\vector(3,-4){40}}               
\put(40, 100){\vector(3,-4){40}}                
\put(120, 100){\vector(3,-4){30}}               

\put(75,70){$\matrix{&&      &\cr
                     &{\cal S} &\sgpd &{\cal V}\cr
                     &&&\cr
                     &\svgpd          &&\svgpd          \cr
                     &&&\cr
                     &{\cal H} &\sgpd & P\cr}$}

\put(100,0){(a)}


\put(200,150){$\matrix{&&      &\cr
                      &T{\cal S}&\sgpd &T{\cal V}\cr
                      &&&\cr
                      &\svgpd          & &\svgpd          \cr
                      &&&\cr
                      &T{\cal H} &\sgpd &TP \cr}$}

\put(240, 160){\vector(3,-4){40}}                
\put(310, 160){\vector(3,-4){40}}               
\put(240, 100){\vector(3,-4){40}}                
\put(320, 100){\vector(3,-4){30}}               

\put(275,70){$\matrix{&&      &\cr
                     &{\cal S} &\sgpd &{\cal V}\cr
                     &&&\cr
                     &\svgpd          &&\svgpd          \cr
                     &&&\cr
                     &{\cal H} &\sgpd & P\cr}$}
\put(300,0){(b)}
\end{picture}\caption{\ \label{fig:triples}}
\end{figure}
We call such structures {\em ${\cal LA}$--double groupoids}. The three
double structures which involve $T^*{\cal S}$ or $T{\cal S}$ we call the
{\em upper faces} or {\em upper structures}, the other three being the
{\em lower faces} or {\em lower structures}. That $\pi_{\cal S}^\#$ respects
all twelve ordinary structures follows from the hypothesis and basic results
for ordinary Poisson groupoids.

In what follows we will use the following result repeatedly. The proof is
straightforward.

\begin{prop}                                   \label{prop:repeat}
\begin{enumerate}
\item The core of each upper face in an \LAdgpd\ is a Lie groupoid or Lie
algebroid over the core of the opposite lower face, and forms an \LAgpd\
with respect to the structures of the bases of the two cores.

\item The restriction of a morphism of \LAdgpd s to the core of an upper
face is a morphism of \LAgpd s.
\end{enumerate}
\end{prop}

The map $\tila_{*{{\cal H}}}\colon A^*_V{\cal S}\to T{\cal H}$ is a
morphism of the \LAgpd s which form the bottom faces in
Figure~\ref{fig:triples}, the other maps being
$a_{*{{\cal C}}}\colon A^*{\cal C}\to TP$ and $\id_{{\cal H}}$. Its core map
$A^*{\cal V}\to A{\cal H}$ is therefore a morphism of Lie algebroids: denote
it by $D_{{\cal H}}$.

By \cite[\S4--\S5]{Mackenzie:1992}, the anchor
$\tila_{{\cal V}}\colon A_V{\cal S}\to T{\cal H}$ is a morphism of
\LAgpd s over $a_{{\cal V}}\colon A{\cal V}\to TP$ and $\id_{\cal H}$, with
core morphism $\da_{A{\cal H}}\colon A{\cal C}\to A{\cal H}$. Its dual
$\tila_{{\cal V}}^*\colon T^*{\cal H}\to A^*_V{\cal S}$ is therefore, by
\ref{prop:dualmorph}, a
morphism over $\da^*_{A{\cal H}}\colon A^*{\cal H}\to A^*{\cal C}$ and
$\id_{\cal H}$, with core map $a^*_{{\cal V}}\colon T^*P\to A^*{\cal V}$.
It follows that $\pi^\#_{\cal H} = \tila_{*{{\cal H}}}\circ\tila^*_{{\cal V}}$
is a morphism of groupoids over
$a_{*{{\cal C}}}\circ \da^*_{A{\cal H}}\colon A^*{\cal H}\to TP$; since it is
also a morphism of Lie algebroids over ${\cal H}$, it is a morphism of
\LAgpd s, and the core morphism is
$D_{{\cal H}}\circ a^*_{{\cal V}}\colon T^*P\to A{\cal H}$. This proves
the first part of the following theorem.

\begin{thm}
With the induced structures, ${\cal H}\gpd P$ is a Poisson groupoid with
$$
a_{*{{\cal H}}} = a_{*{{\cal C}}}\circ\da^*_{A{\cal H}} =
         -a_{{\cal V}}\circ D^*_{{\cal H}}.
$$
The induced Poisson structure on $P$ coincides with that induced by ${\cal C}$.
\end{thm}

For the last statement, recall (as part of \ref{thm:cores}, or see
\cite[\S5, p.230]{Mackenzie:1992}) that
$a_{{\cal C}} = a_{{\cal H}}\circ\da_{A{\cal H}}$. Hence
$a_{*{{\cal C}}}\circ a^*_{{\cal C}} =
a_{*{{\cal C}}}\circ\da^*_{A{\cal H}}\circ a^*_{{\cal H}} =
a_{*{{\cal H}}}\circ a^*_{{\cal H}}.$ The second
equation for $a_{*{{\cal H}}}$ follows by noting that the core map for
$\pi^\#_{\cal H}$ is $D_{{\cal H}}\circ a^*_{{\cal V}}$ and the negative dual
of this is equal to the base map, by \ref{ex:rmk31}.

The same process can be carried out with ${\cal H}$ and ${\cal V}$
interchanged. We now have two Lie algebroid morphisms,
$D_{{\cal H}}\colon A^*{\cal V}\to A{\cal H}$ and
$D_{{\cal V}}\colon A^*{\cal H}\to A{\cal V}$, defined as the cores of
$\tila_{*{{\cal H}}}$ and $\tila_{*{{\cal V}}}$ respectively.

\begin{prop}                                   \label{prop:DD}
$D^*_{{\cal V}} = - D_{{\cal H}}.$
\end{prop}

\pf
We can regard $\tila_{*{{\cal V}}}$ as the base map for the cotangent
\LAgpd s which form the top faces of Figure~\ref{fig:triples}. From
\ref{ex:rmk31} it follows that the core map for these faces is
$-\tila^*_{*{{\cal V}}}\colon T^*{\cal V}\to A_{{\cal H}}{\cal S}$. In turn,
as in \ref{prop:repeat}, the cores of the top faces of
Figure~\ref{fig:triples} are groupoids over the cores of the bottom faces.
The base map $A^*{\cal V}\to A{\cal H}$ for $-\tila^*_{*{{\cal V}}}$ is, by
\ref{prop:dualmorph}, the negative dual of the core map of
$\tila_{*{{\cal V}}}$; that is, it is $-D^*_{{\cal V}}$.

By the commutativity properties of the triple structures, one can see that
the base map of the core map of the top faces is the same as the core map
for the bottom faces. But the core map for the bottom faces is $D_{{\cal H}}$.
\boom

For any Lie algebroid $A$, let $\Bar{A}$ denote the same vector bundle
with bracket $[X,Y]^- = -[X,Y]$ and anchor $\Bar{a} = -a$. If now
${\cal B} = (A, A^*)$ is a Lie bialgebroid, define the {\em flip} of
${\cal B}$ to be $(A^*, \Bar{A})$. The flip of ${\cal B}$ induces on
the base the same Poisson structure as does ${\cal B}$.

Proposition \ref{prop:DD} can now be restated:
$D_{\cal V}\colon A^*{\cal H}\to A{\cal V}$ is a morphism
of Lie bialgebroids from $(A^*{\cal H}, \Bar{A{\cal H}})$ to
$(A{\cal V}, A^*{\cal V}).$

We now specialize to symplectic double groupoids, as considered in
\cite{CosteDW}, \cite{Weinstein:1988} and \cite{LuW:1989}. We take the signs
in each structure on $S$ to be the same.

\begin{df}
A {\em symplectic double groupoid} is a double Lie groupoid
$({\cal S};{\cal H},{\cal V};P)$ together with a symplectic structure
on ${\cal S}$ such that both groupoid structures on ${\cal S}$ are
symplectic groupoids.
\end{df}

Since $\pi_{\cal S}^\#\colon T^*{\cal S}\to T{\cal S}$ is an isomorphism of
double groupoids, it follows that the base maps
$\tila_{*{{\cal H}}}\colon A^*_V{\cal S}\to T{\cal H}$ and
$\tila_{*{{\cal V}}}\colon A^*_H{\cal S}\to T{\cal V}$ are \LAgpd\
isomorphisms, and hence their core maps
$D_{{\cal H}}\colon A^*{\cal V}\to A{\cal H}$ and
$D_{{\cal V}}\colon A^*{\cal H}\to A{\cal V}$ are Lie algebroid isomorphisms.
Further, the core map $\pi^\#_{\cal C}\colon T^*{\cal C}\to T{\cal C}$ of
$\pi^\#_{\cal S}$ is an isomorphism, and so ${\cal C}\gpd P$ is a symplectic
groupoid. This proves the following result.

\begin{thm}                                  \label{thm:sideduality}
Let $({\cal S};{\cal H},{\cal V};P)$ be a symplectic double groupoid with
core ${\cal C}\gpd P$. Then ${\cal C}$ is a symplectic groupoid realizing
$P$ and $D_{\cal V}\colon A^*{\cal H}\to A{\cal V}$ is an isomorphism of
Lie bialgebroids from $(A^*{\cal H}, \Bar{A{\cal H}})$ to
$(A{\cal V},A^*{\cal V}).$
\end{thm}

Define Lie bialgebroids $(A,A^*)$ and $(B,B^*)$ on the same base $P$ to
be {\em dual} if $(A,A^*)$ is isomorphic to $(B^*,\Bar{B})$, and define
Poisson groupoids $G$ and $G'$ on the same base to be {\em dual} if their
Lie bialgebroids $(AG, A^*G)$ and $(AG', A^*G')$ are dual
\cite{Weinstein:1988}. Then \ref{thm:sideduality} shows that
$(A{\cal H}, A^*{\cal H})$ and $(A{\cal V}, A^*{\cal V})$ are dual,
and thus ${\cal H}\gpd P$ and ${\cal V}\gpd P$ are Poisson groupoids in
duality.

\begin{ex}\rm
Let $G\gpd P$ be any Poisson groupoid and let $S = G\times\Bar{G}$ have the
double structure $(S;G,P\times\Bar{P};P)$ of \cite[\S4.5]{Weinstein:1988}
or \cite[2.3]{Mackenzie:1992}; thus $S$ is a pair (Poisson) groupoid over
$G$ and over $P\times\Bar{P}$ it is the Cartesian square of $G$.

Then $AH = AG,\ AV = TP$ and $D_V = a_*,\ D_H = -a^*_*.$ When $G$ is
symplectic, $D_H$ is the isomorphism $T^*P\to AG$ of \cite{CosteDW}, and
$(AG, A^*G)$ is dual to $(TP,T^*P)$.
\end{ex}

Returning to a general symplectic double groupoid, the negative dual of
$\tila_{*{{\cal H}}}$ gives an isomorphism of \LAgpd s from
$(T^*{\cal H};A^*{\cal H},{\cal H};P)$ to $A_V{\cal S}$, the side maps being
$D_{{\cal V}}$ and $\id_{\cal H}$. In particular $A_V{\cal S}\gpd A{\cal V}$
may be identified with the standard symplectic groupoid
$T^*{\cal H}\gpd A^*{\cal H}$, thus giving a symplectic realization of
the linearized Poisson structure on $A{\cal V}$.

Now $A(-\tila^*_{*{\cal H}})\colon AT^*{\cal H}\to A^2{\cal S}$ is
automatically an isomorphism of Lie algebroids over $D_{\cal V}$, and is
easily checked to also be an isomorphism over $A{\cal H}$, since both
vertical structures are prolongations. Referring briefly to
\cite[1.18]{Mackenzie:Doubla2}, there is now an isomorphism
$$
A(-\tila^*_{*{\cal H}})\circ (j'_{\cal H})^{-1}\colon
T^*A{\cal H}\to A^2{\cal S}
$$
of double Lie algebroids. Here $T^*A{\cal H}$ has the cotangent Lie
algebroid structure over $A{\cal H}$ induced by the Poisson structure on
$A{\cal H}$, and over $A^*{\cal H}$ it has the structure transported from
the Lie algebroid of $T^*{\cal H}\gpd A^*{\cal H}$ via $j'_{\cal H}\colon
AT^*{\cal H}\to T^*A{\cal H}$ (see \cite[\S7]{MackenzieX:1994}). We have
proved the following.

\begin{thm}
Let $({\cal S};{\cal H},{\cal V};P)$ be a symplectic double groupoid. Then
the double Lie algebroid $A^2{\cal S} \isom A_2{\cal S}$ is canonically
isomorphic to $(T^*(A{\cal H});A{\cal H},A^*{\cal H};P)$, and to
\newline
$(T^*(A{\cal V});A^*{\cal V},A{\cal V};P)$, with the structures induced
from the Poisson groupoid structures of ${\cal H}$ and ${\cal V}$.
\end{thm}

Finally, for any double Lie groupoid $S$, consider the symplectic double
groupoid ${\cal S} = T^*S$ of \ref{thm:dcot}. As a direct corollary of
\ref{thm:sideduality}, we have the following result.

\begin{thm}                                        \label{thm:pairs}
Let $(S;H,V;M)$ be a double Lie groupoid. Then $A^*_VS\gpd A^*C$ and
\newline
$A^*_HS\gpd A^*C$ are Poisson groupoids in duality.
\end{thm}

For the maps $D_{\cal V}$ and $D_{\cal H}$ in this case, see \ref{prop:DRI}.

Theorem \ref{thm:pairs} provides, at least in principle, a wide class of
examples of dual pairs of Poisson groupoids. Note that \cite{BrownM:1992}
gives a construction principle for double groupoids $S$ in terms of $H, V$
and $C$, subject to conditions analogous to local triviality.

It seems likely that, just as all Poisson manifolds which are compatible with
a vector bundle structure are duals of Lie algebroids, so all dual pairs of
Poisson groupoids which are suitably compatible with an additional linear
structure will arise by an abstract version of \ref{thm:pairs}. This will
be investigated elsewhere.

\section{DUALITY FOR ${\cal LA}$-GROUPOIDS}
\label{sect:dflag}

We must now prove Theorem \ref{thm:needed}. This is an unexpectedly
complicated result, and we obtain it as an instance of a general duality
which is clearer and no more difficult than \ref{thm:needed}. At the
same time this duality reveals a great deal about the internal structure
of \LAgpd s and their Lie algebroids, and this will be crucial to
\cite{Mackenzie:DLADLB}.

We begin with an extension to the general context of \LAgpd s of the results
of \cite{MackenzieX:1994} and \cite{MackenzieX:1998} regarding pairings and
dualities for tangent and cotangent bundles of vector bundles.
In the case of double vector bundles, related results for duals
have been given very recently by Konieczna and Urba\'nski \cite{KoniecznaU};
in particular their Theorem~16 corresponds to our
Theorem~\ref{thm:dualduality}.

Consider a double vector bundle
$$
\matrix{&&\tilq_H&&\cr
        &E&\lrah&E^V&\cr
        &&&&\cr
        \tilq_V&\Big\downarrow&&\Big\downarrow&q_V\cr
        &&&&\cr
        &E^H&\lrah&M,&\cr
        &&q_H&&\cr}
$$
as in \cite[1.1]{Mackenzie:1992}, with core $q_K\colon K\to M$. Thus $E$
is a double Lie groupoid in which all groupoid structures are actually
vector bundles. We can therefore apply the dualization process of Pradines
as in \S\ref{sect:datcdg} to either structure on $E$. Taking the dual of
the vertical structure gives a double vector bundle
\begin{equation}                    \label{eq:vdual}
\matrix{&&\tilq_V^{(*)}&&\cr
        &E^{*V}&\lrah&K^*&\cr
        &&&&\cr
\tilq_{V*}&\Big\downarrow&&\Big\downarrow&q_{K*}\cr
        &&&&\cr
        &E^H&\lrah&M,&\cr
        &&q_H&&\cr}
\end{equation}
with core $(E^V)^*\to M$. Here $E^{*V}\to E^H$ is the usual dual of a
vector bundle. The unfamiliar projection $\tilq_V^{(*)}\colon E^{*V}\to K^*,\
\Phi\mapsto\kappa$ is defined by
\begin{equation}                             \label{eq:unfproj}
\langle\kappa, k\rangle = \langle\Phi, \til0^V_X \plusH \Bar k\rangle
\end{equation}
where $k\in K_m,\ \Phi\colon \tilq_V^{-1}(X)\to\R$ and $X\in E_m^H$. The
addition $\plusH$ in $E^{*V}\to K^*$ is defined by the appropriate form
of (\ref{eq:dual2}), namely
\begin{equation}
\langle\Phi\plusH\Phi', \xi\plusH\xi'\rangle =
   \langle\Phi,\xi\rangle + \langle\Phi',\xi'\rangle
\end{equation}
and the zero above $\kappa\in K^*_m$ is $\til0^{(*V)}_\kappa$ defined by
$$
\langle\,\til0^{(*V)}_\kappa, \til0^H_x \plusV \Bar k\rangle =
\langle\kappa, k\rangle
$$
where $x\in E^V_m, k\in K_m.$ The scalar multiplication is defined in a
similar way. The core element $\Bar\psi$ corresponding to
$\psi\in(E^V_m)^*$ is
$$
\langle\Bar\psi, \til0^H_x \plusV \Bar k\rangle =
\langle\psi,x\rangle.
$$

We call (\ref{eq:vdual}) the {\em vertical dual} of $E$. There is also
of course a {\em horizontal dual}
\begin{equation}                    \label{eq:hdual}
\matrix{&&\tilq_{H*}&&\cr
        &E^{*H}&\lrah&E^V&\cr
        &&&&\cr
\tilq_H^{(*)}&\Big\downarrow&&\Big\downarrow&q_V\cr
        &&&&\cr
        &K^*&\lrah&M,&\cr
        &&q_{K*}&&\cr}
\end{equation}
with core $(E^H)^*\to M$, defined by corresponding formulas. The following
result is now somewhat unexpected.

\begin{thm}                           \label{thm:dualduality}
There is a natural (up to sign) duality between $E^{*V}\to K^*$ and
$E^{*H}\to K^*$ given by
\begin{equation}                       \label{eq:3duals}
\langle\Phi, \Psi\rangle = \langle\Psi, \xi\rangle
                            -  \langle\Phi, \xi\rangle
\end{equation}
where $\Phi\in E^{*V},\ \Psi\in E^{*H}$ have
$\tilq_V^{(*)}(\Phi) = \tilq_H^{(*)}(\Psi)$ and $\xi$ is any element
of $E$ with $\tilq_V(\xi) = \tilq_{V*}(\Phi)$ and
$\tilq_H(\xi) = \tilq_{H*}(\Psi).$
\end{thm}

The pairing on the LHS of (\ref{eq:3duals}) is over $K^*$, whereas the
pairings on the RHS are over $E^H$ and $E^V$ respectively.

\bigskip

\pf
Let $\Phi$ and $\Psi$ have the forms $(\Phi; X, \kappa; m)$ and
$(\Psi; \kappa, x; m)$. Then $\xi$ must have the form $(\xi; X, x; m)$.
If $\eta$ also has the form $(\eta; X, x; m)$ then there is a $k\in K_m$
such that $\xi = \eta \plusV (\til0^V_X \plusH \Bar k)$, and so
$$
\langle\Phi,\xi\rangle = \langle\Phi, \eta\rangle + \langle\kappa, k\rangle
$$
by (\ref{eq:unfproj}). By the interchange law \cite[(2)]{Mackenzie:1992}
we also have $\xi = \eta \plusH (\til0^H_x \plusV \Bar k)$ and so
$$
\langle\Psi,\xi\rangle = \langle\Psi, \eta\rangle + \langle\kappa, k\rangle.
$$
Thus (\ref{eq:3duals}) is well defined. To check that it is bilinear
is routine.

Suppose $\Phi$, given as above, is such that $\langle\Phi, \Psi\rangle = 0$
for all $\Psi\in(\tilq_H^{(*)})^{-1}(\kappa).$ Take any $\phi\in (E^H_m)^*$
and consider $\Psi = \til0_\kappa^{(*H)}\plusH\Bar\phi$. Then, taking
$\xi = \til0^V_X$ we find $\langle\Phi, \xi\rangle = 0$ and
$\langle\Psi, \xi\rangle = \langle\phi, X\rangle.$
Thus $\langle\phi, X\rangle = 0$ for all $\phi\in(E^H_m)^*$ and so
$X = 0^H_m$. It follows \cite[1.2]{Mackenzie:1992} that
$$
\Phi = \til0_\kappa^{(*V)} \plusV\Bar\psi
$$
for some $\psi\in(E^V_m)^*$. Now taking any $k\in K_m$ and defining
$\xi = \til0_x^H\plusV \Bar k$, we find that
$$
\langle\Phi, \xi\rangle = \langle\kappa, k\rangle + \langle\psi, x\rangle
\qquad\mbox{and}\qquad
\langle\Psi, \xi\rangle = \langle\kappa, k\rangle.
$$
So $\langle\psi, x \rangle = 0$ for all $x\in E^H_m$, since a suitable
$\Psi$ exists for any given $x$. It follows that $\psi = 0\in (E^V_m)^*$
and so $\Phi$ is indeed the zero element over $\kappa$. Thus the pairing
(\ref{eq:3duals}) is nondegenerate.
\boom

\begin{ex}\rm
Let $E = TA$ be the tangent double vector bundle of an ordinary vector
bundle $(A, q, M)$ as in \cite[\S5]{MackenzieX:1994}. Regarding $TA\to A$
as the vertical structure we have, in the notation of \cite{MackenzieX:1994},
double vector bundles
$$
\matrix{&&      &&\cr
E^{*V} =&T^*A&\lrah&A^*&\cr
        &&&&\cr
        &\Bigg\downarrow&&\Bigg\downarrow&\cr
        &&&&\cr
        &A&\lrah&M&\cr
        &&&&\cr}
\mbox{\qquad and\qquad}
\matrix{&&      &&\cr
E^{*H} =&T^\sol A&\lrah&TM&\cr
        &&&&\cr
        &\Bigg\downarrow&&\Bigg\downarrow&\cr
        &&&&\cr
        &A^*&\lrah&M&\cr
        &&&&\cr}
$$
Recall also from \cite[\S5]{MackenzieX:1994} the canonical isomorphisms
$R\colon T^*(A^*)\to T^*(A)$ and $I\colon T(A^*)\to T^\sol A$. This $R$
is an isomorphism of double vector bundles over $A$ and $A^*$, and an
antisymplectomorphism, and $I$ is an isomorphism of double vector bundles
over $A$ and $TM$. Using $R$ and $I$ the pairing of $T^*(A)$ and $T^\sol A$
which $E$ induces can be transported to a pairing of $T^*(A^*)$ and $T(A^*)$
and is then
$$
\langle{\cal F}, {\cal X}\rangle = \langle I({\cal X}), \xi\rangle
                                    - \langle R({\cal F}), \xi\rangle,
$$
where ${\cal F}\in T^*(A^*),\ {\cal X}\in T(A^*)$, and $\xi\in TA$ has a
suitable form. Now $\langle I({\cal X}), \xi\rangle$ is precisely the
tangent pairing $\llangle{\cal X}, \xi\rrangle$ of $T(A^*)$ and $TA$ over
$TM$, as in \cite[5.3]{MackenzieX:1994}. And the definition of $R$ in
\cite[(38)]{MackenzieX:1994} is
\begin{equation}                 \label{eq:R}
\langle{\cal F}, {\cal X}\rangle + \langle R({\cal F}), \xi\rangle =
 \llangle {\cal X}, \xi\rrangle
\end{equation}
where the first term is the standard pairing of $T^*(A^*)$ and $T(A^*)$.
It follows that the pairing of $T^*(A^*)$ and $T(A^*)$ induced by $E$ is
precisely the standard one.
\end{ex}

\begin{ex}\rm
Continuing the preceding example, now take $E = T^*A$, the double cotangent
bundle. The two duals are
$$
\matrix{&&      &&\cr
E^{*V} =&TA&\lrah&TM&\cr
        &&&&\cr
        &\Bigg\downarrow&&\Bigg\downarrow&\cr
        &&&&\cr
        &A&\lrah&M&\cr
        &&&&\cr}
\mbox{\qquad and\qquad}
\matrix{&&      &&\cr
E^{*H} =&(T^*A)^\sol&\lrah&A^*&\cr
        &&&&\cr
        &\Bigg\downarrow&&\Bigg\downarrow&\cr
        &&&&\cr
        &TM&\lrah&M&\cr
        &&&&\cr}
$$
where $(T^*A)^\sol$ denotes the dual taken with respect to the structure
over $A^*$. To represent this in more familiar terms, take the dual of
$R^{-1}\colon T^*A\to T^*(A^*)$ with respect to the structures over $A^*$.
Recalling that $R$ induces $-\id\colon T^*M\to T^*M$ on the cores, and
using \ref{prop:dualmorph}, it follows that
$(R^{-1})^*\colon T(A^*)\to (T^*A)^\sol$ sends elements of the form
$({\cal X};\phi,-x;m)$ to $((R^{-1})^*({\cal X});\phi,x;m)$.

Transporting the pairing induced by $E = T^*A$ to a pairing of
$TA\to TM$ and $T(A^*)\to TM$, we take elements $(\xi; X, x; m)$ of $TA$
and $({\cal X}; \phi, -x; m)$ of $T(A^*)$ and get, for
$(\Phi; X, \phi; m)$ in $T^*A$,
$$
\langle\xi, {\cal X}\rangle = \langle (R^{-1})^*({\cal X}), \Phi\rangle
                                    - \langle\xi, \Phi\rangle
= \langle R^{-1}(\Phi), {\cal X}\rangle - \langle\Phi, \xi\rangle
= \llangle{\cal X}, -\xi\rrangle - \langle\Phi, -\xi\rangle
                                       - \langle\Phi, \xi\rangle,
$$
where we applied (\ref{eq:R}) with $\xi$ replaced by its usual negative
$(-\xi; X, -x; m)$. Thus the pairing induced by $E$ between $TA$ and
$T(A^*)$ over $TM$ is
$$
\langle\xi, {\cal X}\rangle = \llangle{\cal X}, -\xi\rrangle,
$$
the sign being an essential feature.
\end{ex}

\begin{ex}\rm
Let $A, B$ and $K$ be any three vector bundles on $M$. The manifold
$E = A\oplus B \oplus K$ may be regarded as the pullback of
$A\oplus K$ to $B$ and as the pullback of $B\oplus K$ to $A$, and with
these structures is a double vector bundle over $M$ with core $K$. Let
$\Phi = X\oplus\psi\oplus\kappa$ be an element of
$E^{*V} = A\oplus B^*\oplus K^*$ and let $\Psi = \phi\oplus x\oplus\kappa$
be an element of $E^{*H}.$ Then taking any $\xi = X\oplus x\oplus k\in E$,
we find that
$$
\langle \Phi, \Psi\rangle = \langle\phi, X\rangle - \langle\psi, x\rangle.
$$
\end{ex}

Although we have proved that $E^{*V}$ and $E^{*H}$ are dual as vector
bundles over $K^*$, we have not yet considered the relationships between the
other structures present. This is taken care of by the following result,
whose proof is straightforward.

\begin{prop}                                         \label{prop:dudvb}
Let $(D;A,D^V;M)$ and $(E;A,E^V;M)$ be double vector bundles with
$D^H = E^H = A$ and with cores $K$ and $L$ respectively. Suppose given a
nondegenerate pairing $\langle\ ,\ \rangle$ of the vertical bundles $D\to A$
and $E\to A$ such that
\begin{enumerate}
\item for all $x\in D^V,\ \ell\in L$, $\langle \til0^H_x,\Bar{\ell}\rangle
= \langle x, \ell\rangle$, a nondegenerate pairing of $D^V$ and $L$ over $M$;
\item for all $k\in K,\ y\in E^V$, $\langle \Bar{k}, \til0^H_y\rangle
= \langle k, y\rangle$, a nondegenerate pairing of $K$ and $E^V$ over $M$;
\item for all $k\in K,\ \ell\in L$, $\langle \Bar{k}, \Bar{\ell}\rangle = 0;$
\item for all $d_1, d_2 \in D,\ \xi_1, \xi_2 \in E$ such that
$\tilq^D_H(d_1) = \tilq^D_H(d_2),\ \tilq^E_H(\xi_1) = \tilq^E_H(\xi_2),\
\tilq^D_V(d_1) = \tilq^E_V(\xi_1),\ \tilq^D_V(d_2) = \tilq^E_V(\xi_2),$
we have $\langle d_1 \plusH d_2, \xi_1 \plusH\xi_2\rangle =
\langle d_1, \xi_1 \rangle + \langle d_2, \xi_2 \rangle$;
\item for all $d\in D, \xi\in E$ such that $\tilq^D_V(d) = \tilq^E_V(\xi)$
and all $t\in\R$, we have
$\langle t\timesH d, t\timesH \xi\rangle = t\langle d, \xi \rangle$.
\end{enumerate}
(In all the above conditions we assume the various elements lie in
compatible fibres over $M$.)

Then the map $F\colon D\to E^{*V},\ (F(d))(\xi) = \langle d, \xi\rangle$
is an isomorphism of double vector bundles, with respect to $\id\colon A\to A$
and the isomorphisms $D^V\to L^*$ and $K\to (E^V)^*$ induced by the pairings
in {\em (i)} and {\em (ii)}.
\end{prop}

Applying this result to the pairing (\ref{eq:3duals}) of $E^{*V}$ and
$E^{*H}$, we find that the induced pairing of $E^H$ and $(E^H)^*$ is the
standard one, but that of $(E^V)^*$ and $E^V$ is the negative of the
standard pairing. Hence there is an unavoidable sign in the following
result.

\begin{cor}                                \label{cor:56}
The pairing {\em (\ref{eq:3duals})} induces an isomorphism of double vector
bundles from \newline
$(E^{*V}; K^*, E^H; M)$ to $((E^{*H})^{*V};K^*, E^H;M)$ which is
the identity on the side bundles $K^*$ and $E^H$, but $-\id$ on the cores
$(E^V)^*$.
\end{cor}

Now consider a \VBgpd\ $(\Om;G,A;M)$ as in (\ref{eq:VBgpd}), and the
double vector bundle
\newline
$(A\Om;AG,A;M)$ obtained by applying the Lie functor
to the two groupoid structures (see \cite[\S1]{Mackenzie:Doubla2}). All
the preceding results may of course be applied to $A\Om$, but in this
case more is true. Denote the vertical dual of $(A\Om;AG,A;M)$ by
$(A^\sol\Om; AG, K^*; M)$.

\begin{prop}
The double vector bundles
$(A^\sol\Om; AG, K^*; M)$ and $(A(\Om^*);AG,K^*;M)$ are isomorphic, where
$A(\Om^*)$ is the \VBgpd\ associated to $(\Om^*;G,K^*;M)$.
\end{prop}

\pf
The canonical pairing between $\Om$ and $\Om^*$ over $G$ is a Lie groupoid
morphism $\Om \pback{G}\Om^*\to\R$ where
$\Om\pback{G}\Om^*\gpd A\pback{M} K^*$ is the pullback groupoid of $\tilq$
and $\tilq_*$. Applying the Lie functor we get a pairing
$$
\llangle\ ,\ \rrangle =
A(\langle\ ,\ \rangle)\colon A\Om\pback{AG} A(\Om^*)\to\R
$$
which satisfies (iv) and (v) of \ref{prop:dudvb} because
$\llangle\ ,\ \rrangle$ is itself a vector bundle morphism. To prove (i) of
\ref{prop:dudvb}, take $X\in A_m$ and $\phi\in A^*_m$ and write
$$
\til0_X = \left.\frac{d}{dt}\tilone_X\right|_0,\qquad\qquad
\Bar{\phi} = \left.\frac{d}{dt}(t\phi)\right|_0
$$
where the $t\phi$ is itself a core element in $\Om^*$. Then
$$
\llangle\til0_X,\Bar{\phi}\rrangle =
\left.\frac{d}{dt}\langle\tilone_X,t\phi\rangle\right|_0 =
\langle\phi,X\rangle
$$
which proves (i). Likewise, $\llangle\Bar{k}, \til0^{(*)}_\kappa\rrangle
= \langle\kappa,k\rangle$ for $k\in K_m, \kappa\in K^*_m$. The proof of
(iii) is similar.

To show that $\llangle\ ,\ \rrangle$ is nondegenerate, take
$(\Xi;x,X;m)$ in $A\Om$ and assume that $\llangle{\cal X},\Xi \rrangle = 0$
for all $({\cal X};x,\kappa;m)$ in $A(\Om^*)$. We first prove that $X = 0$.
Take any $\phi\in A^*_m$ and define
$$
{\cal X} = \left.\frac{d}{dt}(\til0^{(*)}_{g_t}(t\Bar{\phi})^{-1})\right|_0
$$
where $x = \left.\frac{d}{dt}g_t\right|_0$ and $\til0^{(*)}_{g_t}$ is the
zero of $\Om^*$ above $g_t$. Using (\ref{eq:dual2}) and (\ref{eq:dual3})
we get
$$
\langle\til0^{(*)}_{g_t}(t\Bar{\phi})^{-1},\xi_t\tilone_X\rangle =
\langle\til0^{(*)}_{g_t},\xi_t\rangle +
\langle(t\Bar{\phi})^{-1},\tilone_X\rangle =
0 - t\langle\phi, X\rangle
$$
where $\Xi = \left.\frac{d}{dt}\xi_t\right|_0$. Hence
$\llangle{\cal X},\Xi \rrangle = -\langle\phi, X\rangle$. Since $\phi$
was arbitrary, we have $X = 0$ and so $\Xi = A(\til0)(x) \plusH \Bar{k}$
for some $k\in K_m$. Now take any $\kappa\in K^*_m$ and some
$({\cal X};x,\kappa;m)$ in $A(\Om^*)$. Expanding
$\llangle A(\til0)(x)\plusH\Bar{k},{\cal X}\plusH\til0_\kappa \rrangle$
by (iv) of \ref{prop:dudvb} we get $\llangle{\cal X},\Xi \rrangle =
\langle\kappa, k\rangle$ and so $k = 0$. Thus $\Xi$ is the zero of
$A\Om\to AG$ over $x$.
\boom

This of course extends the pairing of $ATG$ and $AT^*G$ given in
\cite[\S5]{MackenzieX:1998}. Likewise extending the notation of
\cite[5.3]{MackenzieX:1994} we define $I_\Om\colon A(\Om^*)\to A^\sol\Om$
by $\langle I_\Om({\cal X}), \Xi\rangle = \llangle{\cal X}, \Xi\rrangle$,
where ${\cal X}\in A(\Om^*),\ \Xi\in A\Om$.

Now consider $A(\Om^*)$ together with the horizontal dual of $A\Om$,
namely $(A^*\Om; K^*, A; M)$.
We define a pairing denoted $\ddagger\ ,\ \ddagger$ between $A(\Om^*)$
and $A^*\Om$ over base $K^*$ by
\begin{equation}                                \label{eq:twiddle}
\ddagger\!{\cal X}, \Psi\ddagger =
\llangle{\cal X}, \Xi\rrangle -
\langle\Psi, \Xi\rangle_A.
\end{equation}
Here ${\cal X}\in A(\Om^*)$ is to have the form $({\cal X}; X,\kappa;m)$ and
$\Psi\in A^*\Om$ is to have the form $(\Psi;\kappa,x;m)$; the element
$\Xi\in A\Om$ must then have the form $(\Xi;X,x;m)$. The pairing
$\langle\ ,\ \rangle_A$ is the standard duality between $A^*\Om$ and $A\Om$.
Using the method of \ref{thm:dualduality} and the properties of
$\llangle\ ,\ \rrangle$, one proves that $\ddagger\ ,\ \ddagger$ is
nondegenerate and that for $\phi\in A^*, x\in A, \psi\in A^*G, X\in AG,$
$$
\ddagger\Bar{\phi},\til0^{(*)}_x\ddagger = \langle\phi,x\rangle
\qquad\mbox{and}\qquad
\ddagger\!A(\til0^{(*)})(X),\Bar{\psi}\ddagger =
- \langle\psi,X\rangle.
$$
Now define $R = R^{gpd}_\Om\colon A^*\Om^*\to A^*\Om$ by
$$
\ddagger{\cal X}, R({\cal F})\ddagger =
\langle{\cal X}, {\cal F}\rangle,
$$
where ${\cal X}\in A\Om^*,\ {\cal F}\in A^*\Om^*$ and where the RHS is the
standard pairing over $K^*$. Thus we have $\langle{\cal X}, {\cal F}\rangle =
\llangle {\cal X}, \Xi\rrangle - \langle R({\cal F}), \Xi\rangle_A,$
similarly to (\ref{eq:R}). The first part of the next result now follows
from using \ref{prop:dudvb} (despite the apparent inconsistency of signs).

\begin{prop}
\begin{enumerate}
\item The map $R^{gpd}_\Om\colon A^*\Om^*\to A^*\Om$ is an isomorphism of
double vector bundles over $K^*$ and $A$ with core map
$-\id\colon A^*G\to A^*G$.
\item The map $R^{vb}_\Om\colon T^*\Om^*\to T^*\Om$, defined as in
{\em \cite[5.5]{MackenzieX:1994}} in terms of the vector bundle structure
of $\Om$, is a groupoid isomorphism over $R^{gpd}_\Om$.
\end{enumerate}
\end{prop}

As in \cite[5.5]{MackenzieX:1994}, one may prove that $R^{gpd}$ is
anti--Poisson with respect to the duals of the Lie algebroid structures.
$R^{gpd}$ is also given by
\begin{equation}                                   \label{eq:RIduals}
R^{gpd} = (I^{\ddagger}\circ\eps)^{-1},
\end{equation}
where $I^\ddagger$ is the dual of $I_\Om$ over $K^*$ and
$\eps\colon A^*\Om\to (A^\sol\Om)^\ddagger$ is induced, by the method of
\ref{cor:56}, using the opposite of the pairing (\ref{eq:3duals}) applied
to $(A\Om;AG,A;M)$. This shows that, in a sense, $R^{gpd}_\Om$ and
(the inverse of) $I_\Om$ are duals.

We can now describe the duality in \ref{thm:pairs}. Consider a double Lie
groupoid $(S;H,V;M)$. Applying the Lie functor to each of the \LAgpd s
in Figure~\ref{fig:LAgpd}, we obtain Lie algebroids
$A^2S = A(A_VS)\to AV$ and $A_2S = A(A_HS)\to AH$, and both $A^2S$ and
$A_2S$ have double vector bundle structures with sides $AH$ and $AV$ and
core $AC$. The canonical involution in the iterated tangent bundle
$T^2S$ restricts to an isomorphism $j = j_S\colon A^2S\to A_2S$ which
preserves $AH, AV$ and $AC$. See \cite[\S2]{Mackenzie:Doubla2}.

We begin by giving simple alternative formulas for the prolonged pairings.
Define
\begin{equation}
{j'}^V = I_V^{-1}\circ j^{*V}\colon A^*(A_HS)\to A(A^*_VS),\qquad
{j'}^H = j^{*H}\circ I_H\colon A(A^*_HS)\to A^*(A_VS),
\end{equation}
where for example $I_V$ refers to the \VBgpd\ $A_VS$ and $j^{*V}$ is the
dual of $j$ over $AH$. (Although the roles of $H$ and $V$ in a general
double groupoid are interchangeable, the definition of $j$ imposes a lack
of symmetry.) Then the following are easily proved:
\begin{equation}
\llangle{\cal X},\Xi\rrangle_{AH} =
\langle({j'}^V)^{-1}({\cal X}), j(\Xi)\rangle,
\qquad
\llangle{\cal Y},\Phi\rrangle_{AV} =
\langle{j'}^H({\cal Y}), j^{-1}(\Phi)\rangle,
\end{equation}
where ${\cal X}\in A(A^*_VS),\ \Xi\in A^2S$ are in the same fibre over $AH$
and ${\cal Y}\in A(A^*_HS),\ \Phi\in A_2S$ are in the same fibre over $AV$.

Using similar techniques and (\ref{eq:RIduals}), we can now prove the
following, which extends Theorem~7.3 of \cite{MackenzieX:1994}.
Denote $R^{gpd}_{A_HS}$ by $R_H$.

\begin{prop}                              \label{prop:DRI}
Let $(S;H,V;M)$ be a double Lie groupoid. Then the duality isomorphisms
${\cal D}_H = D_{A^*_VS}\colon A^*(A^*_HS)\to A(A^*_VS)$ and
${\cal D}_V = D_{A^*_HS}\colon A^*(A^*_VS)\to A(A^*_HS)$ obtained by
applying Theorem~{\em \ref{thm:sideduality}} to ${\cal S} = T^*S$, are
given by
$$
{\cal D}_H = {j'}^V\circ R_H,\qquad
{\cal D}_V = ({j'}^H)^{-1}\circ R_V.
$$
\end{prop}

Using these isomorphisms, all the usual apparatus available for the Lie
algebroid of a Lie groupoid can be applied to the duals.

In the case $S = M^4$, the maps ${j'}^H$ and $({j'}^V)^{-1}$ are the map
$\alpha\colon TT^*M\to T^*TM$ of Tulczyjew \cite{Tulczyjew}, and
${\cal D}_H$ and ${\cal D}_V$ are the map induced by the standard
symplectic structure on $T^*M$.

We can finally show that the duality between Lie algebroids and vector
bundles with a linear Poisson structure \cite{Courant:1990} extends to a
duality between \LAgpd s and structures which we will call \PVBgpd s.

\begin{df}
A \PVBgpd\ is a \VBgpd\ $(\Ups;G,E;M)$ together with a Poisson structure
on $\Ups$ with respect to which $\Ups\gpd E$ is a Poisson groupoid and
$\Ups\to G$ is a Poisson vector bundle.
\end{df}

Denote the core of $\Ups$ by $L$. There are now three double structures on
$T^*\Ups$. It is a double vector bundle $(T^*\Ups;\Ups, \Ups^*;G)$, a
\VBgpd\ $(T^*\Ups;\Ups,A^*\Ups;E)$ and, using $R^{vb}_\Ups$ and
$R^{gpd}_\Ups$, a \VBgpd\ $(T^*\Ups;\Ups^*,A^*\Ups;L^*)$. These fit
together into a triple structure as in Figure~\ref{fig:PVBdual}(a).
\begin{figure}[htb]
\begin{picture}(350,200)  
\put(0,150){$\matrix{&&      &\cr
                      &T^*\Ups&\sgpd &A^*\Ups\cr
                      &&&\cr
                      &\Bigg\downarrow   & &\Bigg\downarrow     \cr
                      &&&\cr
                      &\Ups&\sgpd &E\cr}$}

\put(40, 160){\vector(3,-4){40}}                
\put(110, 160){\vector(3,-4){40}}               
\put(40, 100){\vector(3,-4){40}}                
\put(120, 100){\vector(3,-4){30}}               

\put(75,70){$\matrix{&&      &\cr
                     &\Ups^* &\sgpd &L^*\cr
                     &&&\cr
                     &\Bigg\downarrow &&\Bigg\downarrow \cr
                     &&&\cr
                     &G &\sgpd & M\cr}$}

\put(100,0){(a)}


\put(200,150){$\matrix{&&      &\cr
                      &T\Ups&\sgpd &TE\cr
                      &&&\cr
                      &\Bigg\downarrow & &\Bigg\downarrow \cr
                      &&&\cr
                      &\Ups &\sgpd &E \cr}$}

\put(240, 160){\vector(3,-4){40}}                
\put(310, 160){\vector(3,-4){40}}               
\put(240, 100){\vector(3,-4){40}}                
\put(320, 100){\vector(3,-4){30}}               

\put(275,70){$\matrix{&&      &\cr
                     &TG &\sgpd &TM\cr
                     &&&\cr
                     &\Bigg\downarrow &&\Bigg\downarrow \cr
                     &&&\cr
                     &G &\sgpd & M\cr}$}
\put(300,0){(b)}
\end{picture}\caption{\ \label{fig:PVBdual}}
\end{figure}
We may apply \ref{lem:simple} to the top face using the anchors
$\tila_*\colon A^*\Ups\to TE$ and $\rho\colon\Ups^*\to TG$ and obtain a
map $b\colon L^*\to TM$, and it then follows that $\pi^\#$ is a morphism of
the triple structures inducing $\tila_*, \rho$ and $b$ as corner maps,
with identities elsewhere.

As in \S\ref{sect:sdgapdg}, the three upper cores have double structures,
namely $(T^*G,A^*G, G;M)$,\break $(T^*E;E,E^*;M)$ and $(T^*L;L,L^*;M)$. For
the top face, which is nonstandard, we regard the core as $T^*L$, embedded
into $T^*\Ups$ via the composite $R^{vb}_\Ups\circ c \circ R_L^{-1}$ where
$c$ is the core embedding of $T^*\Ups^*$. Each of these double structures
has core $T^*M$.

The core maps for two of the upper faces are $-\tila^*_*$ and $-\rho^*$ and
we define the core of the top faces to be $\pi^\#_L\colon T^*L\to TL$,
giving a Poisson structure on $L$. This is a morphism of double vector
bundles over $\id_L$ and $b$ and the Poisson structure is thus linear. We
denote by $\da_{AG}$ the core map $E^*\to AG$ of the front faces; because
$\pi^\#_E = \tila_*\circ \tila^*$ must be skew--symmetric, the core map of
the right faces must be $-\da_{AG}^*$. Each of these morphisms of core
double structures induces $-b^*\colon T^*M\to L$ on its cores.

We now have Lie algebroid structures on $\Ups^*, E^*$ and $L^*$, with
anchors $\rho,\ a_G\circ\da_{AG}$ and $b$. We claim that with respect to
these structures, $\Ups^*$ is an \LAgpd. Consider first the target map
$\tilbeta_*\colon\Ups^*\to L^*$. It commutes with the anchors since $\rho$
is a groupoid morphism over $b$. Consider $\xi\in\Ga_G\Ups^*$ and the
associated linear map
$\ell_\xi\colon\Ups\to\R,\ u\mapsto\langle\xi(qu),u\rangle$ where
$q = q_\Ups\colon\Ups\to G$. Then, from (\ref{eq:dual1}), there exists
$Y\in\Ga L^*$ such that $\tilbeta_*(\xi(g)) = Y(\beta g)$ for all
$g\in G$ if and only if there exists $Y$ such that, for all $g\in G$
and $\phi\in L_{\beta g}$,
\begin{equation}                                     \label{eq:ell}
\ell_\xi(\phi\til0_g) = \langle Y(\beta g), \phi\rangle.
\end{equation}

\begin{lem}
Equation {\em (\ref{eq:ell})} is equivalent to the condition that, for all
$g\in G, \phi\in L_{\beta g}$,
\begin{equation}                                   \label{eq:ell1}
(\delta\ell_\xi)(\phi\,\til0_g) = (\delta\ell_Y)(\phi)\,\til0_{\xi(g)}.
\end{equation}
\end{lem}

\pf
The RHS refers to the top face of Figure~\ref{fig:PVBdual}(a), where the
groupoid structure is transported from $T^*\Ups^*\gpd A^*\Ups^*$. Applying
$(R_\Ups^{vb})^{-1}$ and recalling how $T^*L$ is embedded in the top face
as its core, this is
\begin{equation}                              \label{eq:ell2}
[\delta\ell_\phi(Y(\beta g)) -
q^*_{L^*}(Y(\beta g), \delta\langle\phi,Y\rangle)]\,\til0_{\xi(g)}
\end{equation}
where we extended $\phi$ to a section of $L$ and used
\cite[6.4]{MackenzieX:1994}. As in \cite[6.4]{MackenzieX:1994}, the
second term in $[\dots]$ is the pullback of $\delta\langle\phi,Y\rangle$
along $q_{L^*}$ to $Y(\beta g)\in L^*.$ Note that $\til0_{\xi(g)}$ and
the groupoid multiplication are now in $T^*\Ups\gpd A^*\Ups.$

Now, for the cotangent \VBgpd\ of any Lie groupoid $\Om\gpd A$, an
expression $\omega\,\til0_\xi$ for $\omega\in T^*A,\ \xi \in\Om,$ is the
pullback along $\beta$ of $\omega$ to $\xi$. Hence (\ref{eq:ell2}) is
the pullback
$$
\beta^*_{\Ups^*}(\delta\ell_\phi(Y(\beta g))) -
q^*_{\Ups^*}(\beta^*_G(\delta\langle\phi, Y\rangle)),
$$
all at $\xi(g)$, and condition (\ref{eq:ell1}) is that this equals, for
all $g$ and $\phi$,
\begin{equation}                                         \label{eq:ell3}
(R^{vb}_\Ups)^{-1}(\delta\ell_\xi(\phi\til0_g)) =
\delta\ell_\upsilon(\xi(g)) -
q^*_{\Ups^*}(\xi(g), \delta\langle\upsilon,\xi\rangle)
\end{equation}
where we regard $\phi\til0_g$ as the value at $g$ of the pullback
section $\upsilon\in\Ga_G\Ups$ given by
$\upsilon(g') = \phi(\beta g')\til0_{g'}$. Now using (\ref{eq:ell}), we
have $\langle\upsilon,\xi\rangle(g) = \langle Y, \phi\rangle(\beta g),$
and so the two second terms of (\ref{eq:ell2}) and (\ref{eq:ell3})
are equal. Likewise the two first terms are equal.

The converse is proved using the same techniques.
\boom

Any $\xi\in\Ga_G\Ups^*$ induces the vector field
$H_\xi = \pi^\#(\delta\ell_\xi)$ on $\Ups$ (see \cite[\S6]{MackenzieX:1998}
for example). Since $\pi^\#$ is a morphism of the triple structures, a $\xi$
which satisfies (\ref{eq:ell1}) will induce $H_\xi$ such that
\begin{equation}                                      \label{eq:T}
H_\xi(\phi\til0_g) = H_Y(\phi)\til0_{x(g)}
\end{equation}
where $x$ is the anchor vector field on $G$ for $\xi$.

Now if $\xi_1$ is a second such section, with projection $Y_1$, we have
$$
\ell_{[\xi, \xi_1]} = \{\ell_{\xi}, \ell_{\xi_1}\} =
\langle\delta\ell_{\xi_1}, H_{\xi}\rangle
$$
and it follows from the diagram and (\ref{eq:dualcomp}), (\ref{eq:ell1}),
(\ref{eq:T}) that
$$
\ell_{[\xi, \xi_1]}(\phi\til0_g)
= \langle\delta\ell_{Y_1}(\phi)\til0_{\xi_1(g)},H_Y(\phi)\til0_{x(g)}\rangle
= \langle\delta\ell_{Y_1}(\phi),H_Y(\phi)\rangle + 0
= \langle[Y, Y_1](\beta g),\phi\rangle
$$
as desired. This proves that the target of $\Ups^*$ is a Lie algebroid
morphism \cite{HigginsM:1990a}. The proof for the source is similar.

To prove that the groupoid multiplication $\Ups^* * \Ups^*\to\Ups^*$ is a Lie
algebroid morphism, it is sufficient to prove that the bracket of morphic
sections is morphic.

\begin{lem}
$\xi\in\Ga_G\Ups^*$ is a morphism over some $X\in\Ga L^*$ if and only if
$\ell_\xi\colon\Ups\to\R$ is a morphism; that is, if and only if
\begin{equation}                            \label{eq:ellmor}
\langle\xi(hg),vu\rangle = \langle\xi(h),v\rangle + \langle\xi(g),u\rangle
\end{equation}
for all compatible $v,u\in\Ups$, where
$h = q_{\Ups^*}(v), g = q_{\Ups^*}(u).$
\end{lem}

\pf
The direct statement is immediate. Assume that (\ref{eq:ellmor}) holds.
For any $m\in M$, we know that $\xi(1_m) = \tilone_X + k$ for some
$X\in L^*_m, k\in E^*_m$. Taking any $e\in E_m$, put $v = u = \tilone_e$.
Then, using (\ref{eq:dual3}),
$$
\langle e, k\rangle = \langle\tilone_X + k, \tilone_e\rangle = 0
$$
so $k = 0$. We can therefore define $X\in\Ga L^*$ by $\xi(1_m) =
\tilone_{X(m)}$.

Next put $v = \Bar{\phi}$ and $u = \til0_g$ where $\phi\in L_m$ and
$\beta g  = m$. Using (\ref{eq:dual4}) we have
$$
\langle\xi(g),\Bar{\phi}\,\til0_g\rangle =
\langle\tilone_{X(m)},\Bar{\phi}\rangle + 0 =
\langle\phi,X(m)\rangle
$$
and this shows that $\beta_{\Ups^*}(\xi(g)) = X(\beta g).$ Similarly,
$\xi$ projects to $X$ under $\alpha_{\Ups^*}.$ That $\xi$ preserves the
multiplications follows quickly.
\boom

For $\xi\in\Ga_G\Ups^*$ a morphism over $X\in\Ga L^*$ we can now apply
the Lie functor to $\ell_\xi$ and get $A(\ell_\xi)\colon A\Ups\to\R$
which in turn induces a section of $A^*\Ups\to E$ which we denote
$\xi_0$. (Compare \cite[\S5]{MackenzieX:1998}.) The proof of the following
lemma is straightforward.

\begin{lem}
If $\xi\in\Ga_G\Ups^*$ is a morphism over $X\in\Ga L^*$, the 1--form
$\delta\ell_\xi\colon\Ups\to T^*\Ups$ is a morphism over $\xi_0$.
\end{lem}

It follows that $H_\xi = \pi^\#(\delta\ell_\xi)$ is a morphic vector
field on $\Ups$ over $\tila_*\circ\xi$ on $E$. If now $\eta$ is another
morphic section of $\Ups^*$, over $Y\in\Ga L^*$, it follows that
$\ell_{[\xi, \eta]} = \langle\delta\ell_{\eta}, H_{\xi}\rangle\colon
\Ups\to\R$ is a morphism.

This completes the proof of the first half of the following result.
The converse is proved using an adaptation of the above techniques.

\begin{thm}                            \label{thm:LAPVB}
Given a \PVBgpd\ $(\Ups;G,E;M)$ with core $L$ as above, the dual \VBgpd\
$(\Ups^*;G,L^*;M)$ is an \LAgpd\ with core $E^*$.

Conversely, if $(\Om;G,A;M)$ is an \LAgpd\ with core $K$, the dual \VBgpd\
\break $(\Om^*;G,K^*;M)$ is a \PVBgpd\ with core $A^*$.
\end{thm}

\begin{ex}\rm
For any Poisson groupoid $G\gpd P$, the tangent groupoid $(TG;G,TP;P)$ is
a \PVBgpd\ and induces on its core $AG$ the Poisson structure dual
to $A^*G$. (Use \cite[7.3]{MackenzieX:1994}.) Thus its dual $(T^*G;G,A^*G;P)$
is the \LAgpd\ used in \cite{MackenzieX:1998} and \cite{Mackenzie:Doubla2}.
\end{ex}

We can finally prove Theorem \ref{thm:needed}. Recall from
\cite[4.15]{Xu:1995} that if $\phi\colon G\to G'$ is a morphism of Poisson
groupoids then $A(\phi)\colon AG\to AG'$ is a Poisson map with respect
to the Poisson structures dual to the Lie algebroid structures
on $A^*G$ and $A^*G'$.

It follows that if $({\cal S};{\cal H}, {\cal V};P)$ is a Poisson
double groupoid, then its \LAgpd\
\newline
$(A_V{\cal S}, {\cal H}, A{\cal V};P)$ is a \PVBgpd\ and so, by
\ref{thm:LAPVB}, its dual $(A^*_V{\cal S}; {\cal H}, A^*{\cal C}; P)$ is
an \LAgpd.

\newpage


\newcommand{\noopsort}[1]{} \newcommand{\singleletter}[1]{#1}

\end{document}